\documentclass{article}

\usepackage{amssymb,amsfonts,amsthm, amscd}
\usepackage{amsmath}
\usepackage{hyperref}

\numberwithin{equation}{section}
\newtheorem{theorem}{Theorem}[section]

\theoremstyle{remark}

\theoremstyle{definition}

\begin{document}

\begin{center}
\textbf{\Large{Dirac operator with linear potential \\ and its perturbations}}

\vskip 7pt

Yu.A. Ashrafyan and T.N. Harutyunyan

\vskip 7pt

Department of Mathematics and Mechanics, Yerevan State University, Yerevan, Armenia
\end{center}

\begin{abstract}We prove that canonical Dirac expression with linear potential generates operators on axis and half axis,
for which we can find the eigenvalues and eigenfunctions in explicit form.
We construct the perturbations of these operators with in advance given spectra.
\end{abstract}

\textbf{MSC2010:}  Primary 34A55; Secondary 34B30; 47E05

\section{The operator on whole axis}\label{sec1}

The system of differential equations
\begin{equation}\label{eq1.1}
\quad  \ell y\equiv \Big\{ B \frac{d}{dx} + \Omega(x) \Big\} y =\lambda y,
                                                                        \quad y=\left(
                                                                          \begin{array}{c}
                                                                            y_1 \\
                                                                             y_2 \\
                                                                           \end{array}
                                                                        \right),
\end{equation}
is called Dirac canonical system (see \cite{GaLe}, \cite{Ma1} ), if
\begin{equation}\label{eq1.2}
  B = \left(
                         \begin{array}{cc}
                             0 & 1 \\
                             -1 & 0 \\
                           \end{array}
                            \right)
, \qquad
 \Omega (x) = \left(
                                             \begin{array}{cc}
                                               p(x) & q(x) \\
                                               q(x) & -p(x) \\
                                             \end{array}
                                           \right).
\end{equation}
Matrix-function $\Omega(\cdot)$ is usually called the potential.
We assume that $p, q$ are real-valued, local integrable functions.
$\lambda$ is a complex (spectral) parameter, $\lambda \in \mathbb{C}$.

By $L(p, q)$ we denote a self-adjoint operator (see \cite{Na}), generated by differential expression $\ell$ in Hilbert space of two-component vector-functions $L^2( (- \infty, \infty); {\mathbb{C}}^2)$ on the domain
\[
D =  \Big\{ y=\left(         \begin{array}{c}
                                                                     	y_1 \\
                                                                     	y_2
                                                                     \end{array}
                                                                   \right);
y_k \in L^2 (- \infty, \infty) \cap AC (- \infty, \infty);
\]
\begin{equation}\label{eq1.3}
(\ell y)_k \in L^2 (- \infty, \infty), k = 1, 2 \Big\},
\end{equation}
where $AC (- \infty, \infty) = AC (\mathbb{R})$ is the set of functions, which are absolutely continuous on each finite segment
$[a, b] \subset (-\infty, \infty), -\infty < a < b < \infty$.

At first we consider an operator $L(0, x)$ (with $p(x) \equiv 0 $ and $q(x) \equiv x $), which corresponds to the system

\begin{equation}\label{eq1.4}
\quad  \ell y \equiv \Big\{  B \frac{d}{dx} + \Omega_0(x) \Big\} y = \lambda y,
\end{equation}
where $\Omega_0(x) = \left(
                 \begin{array}{cc}
                   0 & x \\
                   x & 0 \\
                 \end{array}
               \right),
$
on the domain \eqref{eq1.3}.
This operator we call Dirac operator with linear potential.

\textbf{Definition.} The values of the parameter $\lambda$, for which the system \eqref{eq1.1} has non trivial solutions from $D \subset L^2((-\infty, \infty);\mathbb{C}^2)$ are called eigenvalues and the corresponding solutions are called eigenfunctions of the operator $L(p, q)$.

As it follows from the results of \cite{Mar} and \cite{LeSa} the spectra of this operator is pure discrete and consists of simple eigenvalues.
The eigenvalues of the operator $L(p, q)$ we will denote $\lambda_n (p, q)$ with corresponding enumeration.
One of the sufficient conditions for discreteness of the spectra is (see \cite{Mar})

\begin{equation}\label{eq1.5}
  \underset{|x| \rightarrow \infty} \lim \left[ p^2(x) + q^2(x) - \sqrt{(p'(x))^2 + (q'(x))^2} \right] = \infty .
\end{equation}
It is easy to see that in our case $(p(x) \equiv 0, q(x) \equiv x)$ the condition \eqref{eq1.5} holds.

Writing the system \eqref{eq1.4} componentwise:

\begin{gather}
  -y_1' + x y_1 = \lambda y_2, \label{eq1.6}\\
  y_2' + x y_2 = \lambda y_1, \label{eq1.7}
\end{gather}
we can obtain two second order differential equations for both $y_1$ and $y_2$ separately:

\begin{gather}
  - y_1'' + x^2 y_1 = (\lambda^2 - 1) y_1, \label{eq1.8}\\
  - y_2'' + x^2 y_2 = (\lambda^2 + 1) y_2. \label{eq1.9}
\end{gather}
It is well known (see, e.g. \cite{LeSa}) that the equation

\[
  - y '' + x^2 y = \mu y
\]
has solution from $L^2(-\infty, \infty)$ only for $\mu = 2 n + 1, \ n = 0, 1, 2, \ldots,$  and corresponding solutions are
Chebyshev-Hermite polynomials

\begin{equation}\label{eq1.10}
  H_n(x) = (-1)^n e^{x^2} \cfrac{d^n e^{-x^2}}{d x^n}.
\end{equation}
Therefore, $\lambda$ can be eigenvalue of $L(0, x)$ only if $\lambda^2 - 1 = 2 n + 1$, i.e. $\lambda^2 = 2 (n + 1)$.
Thus, if $\lambda = \lambda_{\pm n} = \pm \sqrt{2 (n + 1)}$, then the solutions of the equation \eqref{eq1.8} are

\[
y_1(x) = H_n(x), \qquad n = 0, 1, 2\ldots.
\]

At the same time $\lambda^2 + 1 = 2 n + 3 = 2 (n + 1) +1$ and consequently the solutions of the equation \eqref{eq1.9} are

\[
y_2(x) = H_{n+1}(x), \qquad n = 0, 1, 2\ldots.
\]

The Chebyshev-Hermite polynomials have the properties (see \cite{LeSa})
\[
{H'}_n(x) = 2 n H_{n-1}(x), \quad H_{n+1}(x)-2x H_n(x) + 2n H_{n-1}(x) = 0 \qquad n = 1, 2, \ldots.
\]
The general formulae for $H_n(x)$ are
\[
H_n(x) = (2x)^n - \frac{n(n-1)}{1!} (2x)^{n-2} + \frac{n(n-1)(n-2)(n-3)}{2!} (2x)^{n-4} + \cdots ,
\]
in which the last member is $(-1)^{\frac{n}{2}}\frac{n!}{(n/2)!}$, for even $n$ and
$(-1)^{\frac{(n-1)}{2}}\frac{n!}{((n-1)/2)!} 2x$, for odd $n$.
Thus, we note that $H_{2k+1}(0) = 0$, for $k = 0, 1, 2, \ldots$.

It is well known (see, e.g. \cite{LeSa}) that the squares of the $L^2$-norm of $H_n(x)$ with the weight $e^{-x^2}$ is equal
\[
\int_{-\infty}^{\infty} H_n^2(x) e^{-x^2} dx = 2^n n! \sqrt{\pi}
\]
and
\[
\int_{-\infty}^{\infty} H_n(x) H_m(x) e^{-x^2} dx = 0, \qquad n, m = 0, 1, 2, \ldots, \quad n \neq m.
\]
Therefore, if we take
\begin{equation}\label{eq1.11}
  \varphi_n (x) = C_n e^{- \frac{x^2}{2}} H_n(x), \qquad n = 0, 1, 2, \ldots,
\end{equation}
where
\begin{equation}\label{eq1.12}
  C_n = \frac{1}{\sqrt{2^n n! \sqrt{\pi}}},
\end{equation}
then the system $\{ \varphi_n(x) \}_{n=0}^{\infty}$ will became orthonormal system on whole real axis.
It is called the system of Chebyshev-Hermite orthonormal functions.
Now let us show that vector-functions
\[
  U_{-n}(x) = \left(
         \begin{array}{c}
           - \varphi_{n-1}(x) \\
           \varphi_n(x) \\
         \end{array}
       \right),
\]
\begin{equation}\label{eq1.13}
  U_0(x) = \left(
         \begin{array}{c}
           0 \\
           \varphi_0(x) \\
         \end{array}
       \right),
\end{equation}
\[
  U_{n}(x)  = \left(
         \begin{array}{c}
           \varphi_{n-1}(x) \\
           \varphi_n(x) \\
         \end{array}
       \right),
\]
for $n = 1, 2, \ldots$, corresponding to eigenvalues $\lambda_{-n} = - \sqrt{2n}, \lambda_0 = 0, \lambda_n = \sqrt{2n}$, are eigenfunctions of the operator $L(0, x)$.
At first we will show that for $n = 1, 2, \ldots :$

\begin{equation}\label{eq1.14}
  \Bigg \{ \begin{array}{c}
    \varphi^{\prime}_n(x) + x \varphi_n (x) = \sqrt{2n} \varphi_{n-1} (x), \\
    \\
    - \varphi^{\prime}_{n-1} (x) + x \varphi_{n-1} (x) = \sqrt{2n} \varphi_n (x).
  \end{array}
\end{equation}

Indeed, for $\varphi^{\prime}_n (x), n = 1, 2, \ldots$, we have
\[
\varphi^{\prime}_n (x) = - C_n x e^{- \frac{x^2}{2}} H_n(x) + C_n e^{- \frac{x^2}{2}} H^{\prime}_n(x) =
C_n e^{- \frac{x^2}{2}} \left( H^{\prime}_n (x) - x H_n(x) \right),
\]
Putting these into the left side of the equation \eqref{eq1.7} and using the property ${H'}_n(x) = 2 n H_{n-1}(x)$ we will get equalities

\begin{gather*}
\varphi^{\prime}_n (x) + x \varphi_n (x) = \\
C_n e^{- \frac{x^2}{2}} \left( H^{\prime}_n (x) - x H_n (x) \right) + x C_n e^{- \frac{x^2}{2}} H_n (x) = \\
C_n e^{- \frac{x^2}{2}} H^{\prime}_n (x) =
C_n e^{- \frac{x^2}{2}} 2 n H_{n-1} (x) = \\
\cfrac{C_n}{C_{n-1}} 2 n C_{n-1} e^{- \frac{x^2}{2}} H_{n-1} (x) =
\cfrac{2 n C_n}{C_{n-1}} \varphi_{n-1} (x).
\end{gather*}
Taking into account \eqref{eq1.12}, we see that the fraction $\frac{2 n C_n}{C_{n-1}} = \sqrt{2 n}$.
Thus, we have

\[
  \varphi^{\prime}_n (x) + x \varphi_n (x) = \sqrt{2 n} \varphi_{n-1} (x), \qquad n = 1, 2, \ldots .
\]
In the similar way we obtain the following equations (here we use the property $H_{n+1}(x)-2x H_n(x) + 2n H_{n-1}(x) = 0$ )

\[
  - \varphi^{\prime}_{n-1} (x) + x \varphi_{n-1} (x) = \sqrt{2 n} \varphi_n (x), \qquad n = 1, 2, \ldots.
\]
Thus, we have $\ell U_n (x) = \sqrt{2n} U_n (x), \ n=1, 2, \ldots$, i.e.
$U_n (x), \ n = 1, 2, \ldots ,$ are the eigenfunctions of the operator $L(0, x)$ with the eigenvalues $\lambda_n (0, x) = \sqrt{2 n}, \ n = 1, 2, \ldots$.

$U_{-n}(x)$ will satisfy to the system

\begin{equation}\label{eq1.15}
  \Bigg \{ \begin{array}{c}
    \varphi^{\prime}_n(x) + x \varphi_n (x) = (- \sqrt{2 n}) ( - \varphi_{n-1} (x)), \\
    \\
    - (- \varphi^{\prime}_{n-1} (x)) + x (- \varphi_{n-1} (x)) = (- \sqrt{2 n}) \varphi_n (x).
  \end{array}
\end{equation}
In fact the systems \eqref{eq1.15} and \eqref{eq1.14} coincide, which means that for $n = 1, 2, \ldots$
$U_{-n} (x)$ are also the solutions (eigenfunctions) for the system \eqref{eq1.14} (\eqref{eq1.15})
with the eigenvalues $\lambda_{-n} (0, x) = - \sqrt{2 n}, \ n = 1, 2, \ldots$.

$U_0(x)$ satisfies to the system \eqref{eq1.4}, when $\lambda_0 (0, x) = 0$
(note that $\varphi_0(x) = \cfrac{1}{\pi^{\frac{1}{4}}} e^{- \frac{x^2}{2}}$ and
$\varphi^{\prime}_0(x) = \cfrac{1}{\pi^{\frac{1}{4}}} (-x) e^{- \frac{x^2}{2}}$):

\[
  \Bigg \{ \begin{array}{c}
    \varphi^{\prime}_0(x) + x \varphi_0 (x) = \cfrac{1}{\pi^{\frac{1}{4}}} (-x) e^{- \frac{x^2}{2}} +
    x \cfrac{1}{\pi^{\frac{1}{4}}} e^{- \frac{x^2}{2}} = 0  , \\
    -0 + x \cdot 0 = 0 \cdot \cfrac{1}{\pi^{\frac{1}{4}}} e^{- \frac{x^2}{2}}.
  \end{array}
\]

So, such defined vector-functions $U_n(x), \ n \in \mathbb{Z}$ are eigenfunctions of the operator
$L(0, x)$ with the eigenvalues $\lambda_{n} (0, x) = \sqrt{2 |n|} sign(n), \ n  \in \mathbb{Z} $.

\

\section{The operators on half axis}\label{sec2}
For $\alpha \in \big( -\frac{\pi}{2}, \frac{\pi}{2} \big]$, by $L(p, q, \alpha)$ we denote the self-adjoint operator, generated by differential expression $\ell$ \ (see \eqref{eq1.1}) in Hilbert space
of two component vector-functions $L^2(( 0, \infty); {\mathbb{C}}^2)$ on the domain
\[
D_{\alpha} =  \Big\{ y=\left(                                               \begin{array}{c}
                                                                     	y_1 \\
                                                                     	y_2
                                                                     \end{array}
                                                                   \right);
y_k \in L^2 (0, \infty) \cap AC (0, \infty);
\]
\begin{equation}\label{eq2.1}
\quad \ \quad \ (\ell y)_k \in L^2 ( 0, \infty), k = 1, 2; \ y_1(0) \cos \alpha + y_2(0) \sin \alpha = 0 \Big\},
\end{equation}
where $AC (0, \infty)$ is the set of functions, which are absolutely continuous on each finite segment
$[a, b] \subset (0, \infty), 0 < a < b < \infty$.
The eigenvalues of such an operator we will denote by $\lambda_n (p, q, \alpha)$ (in corresponding enumeration).

It is easy to see that if in boundary condition $y_1(0) \cos \alpha + y_2(0) \sin \alpha = 0$ we take $\alpha = 0$, then we have condition

\begin{equation}\label{eq2.2}
  y_1 (0) = 0,
\end{equation}
and if we take $\alpha = \frac{\pi}{2}$, we obtain boundary condition

\begin{equation}\label{eq2.3}
  y_2 (0) = 0.
\end{equation}

It is easy to see from \eqref{eq1.10}-\eqref{eq1.13} that the eigenfunctions of the operator $L(0, x, 0)$ are
vector-functions $U_{2k}(x)$, which correspond to the eigenvalues $\lambda_{k} (0, x, 0) = \lambda_{2 k} (0, x) = 2 \sqrt{|k|} sign(k), \ k \in \mathbb{Z}$.
And the eigenfunctions of the operator $L(0, x, \frac{\pi}{2})$ are vector-functions $U_{2k+1}(x)$ which correspond to the eigenvalues $\lambda_{k} (0, x, \pi / 2) = \lambda_{2 k + 1}(0, x) = \sqrt{2 |2 k + 1|} sign(2 k + 1) , \ k \in \mathbb{Z}$.

By $\psi = \psi(x, \lambda, \alpha, \Omega)$ we denote the solution of the Cauchy problem
($\alpha \in \mathbb{C}$)

\begin{equation}\label{eq2.4}
  \ell y  =  \lambda y, \qquad
    y(0) = \left(
           \begin{array}{c}
             \sin \alpha \\
             - \cos \alpha \\
           \end{array}
         \right),
\end{equation}
on $(0, \infty)$, and we denote this problem by $S(p, q, \lambda, \alpha)$.
Such solution exists and unique and its components $\psi_1$ and $\psi_2$ are entire functions in parameters $\lambda$ and $\alpha$ (see, e.g. \cite{Ha2}).

If $\alpha = 0$ and $\Omega = \Omega_0$, then $\psi (x, \lambda, 0, \Omega_0)$ satisfies to the boundary condition \eqref{eq2.2} and in order to be an eigenfunction of the operator $L(0, x, 0)$ it must be from $L^2 (0, \infty; \mathbb{C}^2)$.
As we have seen recently, it is possible only when $\lambda = \lambda_{2k}(0, x) = 2 \sqrt{|k|} sign(k), \ k \in \mathbb{Z}$.
Thus the eigenvalues and eigenfunctions of the operator $L(0, x, 0)$ are $\lambda_k (0, x, 0)$ and
$\psi(x, \lambda_k(0, x, 0), 0, \Omega_0) = \psi (x, \lambda_{2k}(0, x), 0, \Omega_0)$, for $ k \in \mathbb{Z}$.

Let us now consider Cauchy problems $S(0, x, \lambda_{n}(0, x, 0), 0)$, for $n \in \mathbb{Z}$.
It is easy to see that the functions
\begin{equation}\label{eq2.5}
  V_n(x) = - \cfrac{U_{2n}(x)}{\varphi_{2n}(0)}, \qquad n \in \mathbb{Z},
\end{equation}
are the solutions of the these Cauchy problems.
At the same time $V_n(x)$ are eigenfunctions of the operator $L(0, x, 0)$ which correspond to the eigenvalues $\lambda_n(0, x, 0)$, for $n \in \mathbb{Z}$.
Since the solution of Cauchy problem is unique, it follows that
\begin{equation}\label{eq2.6}
  V_n(x) \equiv \psi(x, \lambda_n(0, x, 0), 0, \Omega_0), \qquad n \in \mathbb{Z}.
\end{equation}

The squares of the $L^2$-norms of these functions

\[
a_n = a_n (0, x) = \| V_n (\cdot)\|^2 =
\displaystyle \int_{0}^{\infty} | V_{n,1}(x) |^2 + | V_{n,2}(x) |^2 dx
\]
are called norming constants.
Using \eqref{eq1.11}-\eqref{eq1.13} and \eqref{eq2.5} we can easily calculate the values of the norming constants:

\[
  a_0 = \cfrac{\pi ^ {1/2}}{2}, \quad
  a_{-n} = a_n = \cfrac{1}{|\varphi_{2n}(0)|^2} = \cfrac{4^n (n!)^2 \pi ^ {1/2}}{(2n)!}, \qquad n = 1, 2, \ldots.
\]
The norming constants and eigenvalues are called spectral data of the operator $L(0, x, 0)$.

Thus, we have two "model" operators on half axis with pure discrete spectra, for which we know eigenvalues, eigenfunctions and norming constants.
Now we want to construct new operators (with in advance given spectra) on half axis, starting from these "model" operators.

\

\section{On the changing of the spectral function in finitely many points}\label{sec3}

The spectral function of an operator $L(0, x, 0)$ is defined as \cite{GaLe, LeSa}

\[
 \rho(\lambda) = \Bigg \{ \begin{array}{c}
                                   \sum_{0 < \lambda_n \leq \lambda} \ a_n^{-1}, \qquad \lambda > 0,   \\
                                   - \sum_{\lambda < \lambda_n \leq 0} \ a_n^{-1}, \qquad \lambda < 0,
                                   \end{array}
\]
and $\rho(0) = 0$, i.e. $\rho(\lambda)$ is left-continuous, step function with jumps in points $\lambda = \lambda_n$ equals $a_n^{-1}$.

In what follows $\delta(x)$ is Dirac $\delta$-function (see, e.g. \cite{Sch}), $\delta_{ij}$ is Kronecker symbol and $v_{ij}(x) = \int_0^{x} V_{i}^{*} (s) V_{j}  (s) ds$, where by the sign $^{*}$ we denote a transponation of vector functions, i.e. $\psi^{*} (x, \lambda) = (\psi_{1}(x, \lambda) \ \psi_{2}(x, \lambda))$, (note that $v_{ij}(x)$ is a scalar function).

In this paragraph we will answer the question, what will happen with the potential $\Omega_0(x)$ if we change spectral data, i.e., if we add or subtract eigenvalues and change the values of norming constants. 
It was proved (see \cite{Ha3}), that if $\rho(\lambda)$ is a spectral function of some self-adjoint operator $L(p, q, \alpha)$, then a function $\tilde \rho(\lambda)$, which differs from $\rho(\lambda)$ by only for finite number of points and is still remaining left-continuous, increasing, step function, is also spectral.
It means that there exists a self-adjoint canonical Dirac operator $\tilde  L (\tilde p, \tilde q, \alpha)$, for which $\tilde \rho(\lambda)$ is spectral function.

At first, we want to construct a new operator $\tilde L(\tilde p, \tilde q, 0)$,
which has the same spectra as $L(0, x, 0)$ except one eigenvalue.
For instance, if we extract eigenvalue $\lambda_0(0, x, 0) = 0$ we will get the following

\begin{theorem}\label{thm3.1}
Let $\rho (\lambda)$ is a spectral function of the operator $L(0, x, 0)$. Then the function $\tilde{\rho} (\lambda)$, defined by relation

\[
 \tilde \rho (\lambda) = \Bigg \{ \begin{array}{c}
                                  \rho(\lambda), \qquad \lambda \leq \lambda_0,   \\
                                  \rho(\lambda) - a^{-1}_0, \qquad \lambda > \lambda_0,
                                   \end{array}
\]
where $a_0 = \sqrt{\pi}/2$, i.e.

\begin{equation}\label{eq3.1}
  d \tilde{\rho} (\lambda) = d \rho (\lambda) -\frac{1}{a_0} \delta(\lambda - \lambda_0) d \lambda
\end{equation}
is also spectral. Moreover, there exists unique self-adjoint canonical Dirac operator $\tilde{L}$ generated by the differential
expression $\tilde{l} = B \cfrac{d}{dx} + \tilde{\Omega}(x)$ and the boundary condition \eqref{eq2.2}, for which $\tilde{\rho} (\lambda)$
is spectral function.
Wherein, the potential function $\tilde{\Omega}(x)$ is represented by the following formula

\begin{equation}\label{eq3.2}
\tilde{\Omega}(x) =    \left(
                        \begin{array}{cc}
                          0 & x - \displaystyle \frac{e^{-x^2}}{a_0 - \int_{0}^{x} e^{-s^2} ds} \\
                          x - \displaystyle \frac{e^{-x^2}}{a_0 -  \int_{0}^{x} e^{-s^2} ds} & 0 \\
                        \end{array}
                      \right)
\end{equation}
and for the eigenfunctions we obtain the formulae
\begin{equation}\label{eq3.3}
  \tilde V_n(x) = \left(
                  \begin{array}{c}
                    V_{n,1}(x) \\
                    \\
                    V_{n,2}(x) + \displaystyle \frac{e^{-\frac{x^2}{2}} \int_{0}^{x} e^{-\frac{s^2}{2}} V_{n,2}(s) ds}{a_0 - \int_{0}^{x} e^{-s^2} ds} \\
                  \end{array}
                \right),
                \qquad n \in \mathbb{Z} \backslash \{0\}.
\end{equation}

\end{theorem}

\begin{proof}
At first we denote $\tilde{\psi}(x, \lambda) = \psi(x, \lambda, 0, \tilde{\Omega})$ and $\psi(x, \lambda) = \psi(x, \lambda, 0, \Omega_0)$.
It is known (see \cite{GaLe, Ma1, LeSa, Ha, AlHrMy}), that there exists transformation operator $\mathbb{I}+\mathbb{G}$:

\begin{equation}\label{eq3.4}
\tilde{\psi}(x,\lambda) = (\mathbb{I}+\mathbb{G}) \psi (x,\lambda) = \psi (x,\lambda) + \int_0^x G(x, s) \psi(s, \lambda) ds,
\end{equation}
which transforms the solution $\psi(x, \lambda)$ of the Cauchy problem $S(0, x, \lambda, 0)$ to the solutions $\tilde{\psi}(x,\lambda)$ of the Cauchy problem
$S(\tilde{p}, \tilde{q}, \lambda, 0)$.
It is also known (see, e.g. \cite{GaLe,LeSa}), that the kernel $G(x, y)$ satisfies to the Gel'fand-Levitan integral equation:

\begin{equation}\label{eq3.5}
G(x,y)+F(x,y)+ \int_0^x G(x,s)F(s,y)ds = 0,\quad 0 \leq y \leq x < \infty,
\end{equation}
where matrix function $F(x, y)$ is defined by the formula

\begin{equation}\label{eq3.6}
F(x,y)= \int_{-\infty}^{\infty} \psi(x,\lambda) \psi^{*} (y,\lambda) d[\tilde{\rho}(\lambda) -\rho(\lambda)].
\end{equation}
It is also known that the potentials $\tilde{\Omega}(x)$ and $\Omega_0(x)$ are connected by the relation
\begin{equation}\label{eq3.7}
\tilde{\Omega}(x) = \Omega_0(x) + G(x, x) B - B G(x, x).
\end{equation}

From the \eqref{eq1.10}-\eqref{eq1.13} and definition \eqref{eq2.5} it follows, that $V_0^{*}(x) = (0 \ \ e^{-\frac{x^2}{2}})$.
Putting the relation \eqref{eq3.1} into \eqref{eq3.6}, and using \eqref{eq2.6}, for the kernel $F(x, y) = F_0(x, y)$ we obtain:
\[
  F_0(x,y) = -a_0^{-1} \psi(x,\lambda_0) \psi^{*} (y,\lambda_0) = -a_0^{-1} V_0(x) V_0^{*} (y) =
\]
\begin{equation}\label{eq3.8}
                                                               = \left(
                                                                \begin{array}{cc}
                                                                  0 & 0 \\
                                                                  0 & -a_0^{-1} e^{-\frac{(x^2+y^2)}{2}} \\
                                                                \end{array}
                                                              \right).
\end{equation}
After some calculations from the equation \eqref{eq3.5} and formula \eqref{eq3.8} for $G_0(x,y)$ we obtain
\[
  G_0(x,y) = \cfrac{1}{a_0 - \int_{0}^{x} e^{-s^2} ds} V_0(x) V_0^{*} (y) =  \left(
                                                                \begin{array}{cc}
                                                                  0 & 0 \\
                                                                  0 & \cfrac{e^{-\frac{(x^2+y^2)}{2}}}{a_0 - \int_{0}^{x} e^{-s^2} ds} \\
                                                                \end{array}
                                                              \right).
\]
Now taking into account \eqref{eq2.6}, putting $G_0(x,y)$ into the equations \eqref{eq3.4} and \eqref{eq3.7} we can easily obtain \eqref{eq3.2} and \eqref{eq3.3}.
Theorem \eqref{thm3.1} is proved.

\end{proof}

Now we want to subtract any finite number $n$ of eigenvalues.
For this reason we denote by $Z_n$ the arbitrary set of finite $n$ number of integers, in increasing order, $Z_n = \{z_1, z_2, \ldots, z_n\} \subset \mathbb{Z}$
(e.g., if $Z_4 = \{z_1, z_2, z_3, z_4\} = \{ -127, 0 , 32, 1259\}$, for
$\displaystyle \sum_{i=1}^4 s_{z_i} = s_{-127} + s_0 + s_{32} + s_{1259}$).

\begin{theorem}\label{thm3.2}
Let $\rho (\lambda)$ is the spectral function of the operator $L$.
Then the function $\tilde{\rho} (\lambda)$, defined by relation
\begin{equation}\label{eq3.9}
  d \tilde{\rho} (\lambda) = d \rho (\lambda)  - \sum_{k=1}^n a_{z_k}^{-1} \delta(\lambda - \lambda_{z_k}) d \lambda
  \end{equation}
also is spectral. Moreover, there exists unique self-adjoint canonical Dirac operator $\tilde{L}$ generated on half axis by the differential
expression $\tilde{l} = B \cfrac{d}{dx} + \tilde{\Omega}(x)$ and the boundary condition \eqref{eq2.2}, for which $\tilde{\rho} (\lambda)$
is spectral function.
Wherein, the potential function $\tilde{\Omega}(x)$ is

\begin{equation}\label{eq3.10}
\tilde{\Omega}(x) =   \left(
                        \begin{array}{cc}
                          p(x, n) & q(x, n) \\
                          q(x, n) & -p(x, n) \\
                        \end{array}
                      \right),
\end{equation}
where $ p(x, n)$ and  $q(x, n)$ are defined by the following formulae:
\[
\begin{array}{c}
  p(x, n)  =  - \cfrac{1}{\det S(x, n)} \displaystyle \sum_{k=1}^n \displaystyle \sum_{p=1}^2 V_{z_k,(3-p)}(x) \det S_{p}^{(k)}(x, n) ,\\
  \\
  q(x, n)  = x + \cfrac{1}{\det S(x, n)} \displaystyle \sum_{k=1}^n \displaystyle \sum_{p=1}^2 (-1)^{p-1} V_{z_k,p}(x) \det S_{p}^{(k)}(x, n) ,
\end{array}
\]
where $S(x, n)$ is $n \times n$ square matrix $ S(x, n) = \{ \delta_{z_i z_j} - a_{z_j}^{-1} v_{z_i z_j} (x) \}_{i,j=1}^n $
and $S_{p}^{(k)}(x, n)$ are matrices, which are obtained from the matrix $S(x, n)$, when we replace $k$-th column of $S(x, n)$
by $H_p(x, n) = \{ a_{z_i}^{-1} V_{z_i,p}(x) \}_{i = 1}^n$ column, $p = 1, 2$.
And for the eigenfunctions $\tilde V_m(x)$ ($ \ m \in \mathbb{Z} \backslash Z_n$) we obtain the representations
\[
  \tilde V_m(x) = \left(
                  \begin{array}{c}
                    V_{m,1}(x) + \cfrac{1}{\det S(x,n)} \displaystyle \sum_{k=1}^n v_{z_k m}(x) \det S_{1}^{(k)}(x, n) \\
                    \\
                    V_{m,2}(x) + \cfrac{1}{\det S(x,n)} \displaystyle \sum_{k=1}^n v_{z_k m}(x) \det S_{2}^{(k)}(x, n) \\
                  \end{array}
                \right).
\]

\end{theorem}

\begin{proof}

In this case the kernel $F(x, y)$ can be written in the following form:

\begin{equation}\label{eq3.11}
F(x, y) = F_n(x, y)= \displaystyle \sum_{k=1}^n -a_{z_k}^{-1} V_{z_k}(x) V_{z_k}^{*} (y),
\end{equation}
and consequently, the integral equation \eqref{eq3.5} becomes to an integral equation with degenerate kernel, i.e. it becomes to a system of linear equations and we will look for the solution in the following form:

\begin{equation}\label{eq3.12}
G_n(x, y)= \displaystyle \sum_{k=1}^n g_{z_k}(x) V_{z_k}^{*} (y),
\end{equation}
where $g_{z_k}(x) = \left(
                  \begin{array}{c}
                    g_{z_k,1} (x) \\
                    g_{z_k,2} (x) \\
                  \end{array}
                \right)
$
is unknown vector-function.
Putting the expressions \eqref{eq3.11} and \eqref{eq3.12} into the integral equation \eqref{eq3.5} we will obtain
a system of algebraic equations for determining the vector-functions $g_{z_k}(x)$:

\begin{equation}\label{eq3.13}
g_{z_k}(x) - \displaystyle \sum_{i = 1}^n a_{z_k}^{-1} v_{z_i z_k}(x) g_{z_i}(x) = a_{z_k}^{-1} V_{z_k}(x), \quad k = 1, 2, \ldots, n.
\end{equation}
It would be better if we consider the equations \eqref{eq3.13} for the vectors $g_{z_k} (x)$
by coordinates $g_{z_k,1}(x)$ and $g_{z_k,2}(x)$ to be systems of scalar linear equations:

\begin{equation*}
g_{z_k,p}(x) - \displaystyle \sum_{i = 1}^n  a_{z_k}^{-1} v_{z_i z_k}(x) g_{z_i,p}(x) = a_{z_k}^{-1} V_{z_k,p}(x), \quad k=1, 2, \ldots, n,  \quad p= 1, 2.
\end{equation*}
The latter systems might be written in matrix form

\[
S(x, n) g_p(x, n) = H_p(x, n), \qquad p=1, 2,
\]
where the column vectors $g_p(x, n) = \{ g_{z_k,p} (x, n)\}_{k=1}^n,\ p = 1, 2$.
The solution of this system can be found in the form (Cramer's rule):

\[
g_{z_k,p}(x, d)
 = \frac{\det S_{p}^{(k)}(x, n)}{\det S(x, n)}, \quad k = 1, 2, \ldots, n, \quad p = 1, 2.
\]

Thus we have obtained for $g_{z_k}(x)$ the following representation:
\[
g_{z_k}(x, n) = \frac{1}{\det S(x, n)}
\left(
  \begin{array}{c}
    \det S_{1}^{(k)}(x, n) \\
    \det S_{2}^{(k)}(x, n) \\
  \end{array}
\right).
\]
Using these $g_{z_k}(x, n)$, from \eqref{eq3.12} we find the function $G_n(x, y)$.
Now taking into account \eqref{eq2.6}, putting $G_n(x,y)$ into the equations \eqref{eq3.7} and \eqref{eq3.4} we obtain the representations for $p(x, n)$, $q(x, n)$ and $\tilde{V}_m(x), \ m \in \mathbb{Z} \backslash \{ Z_n \}$.

Theorem \eqref{thm3.2} is proved.
\end{proof}

 The following theorem says that one can change the values of the finite number norming constants $a_n$ by any positive number $b_n \neq a_n$.
 
\begin{theorem}\label{thm3.3}
Let $\rho (\lambda)$ is the spectral function of the operator $L$.
Then the function $\tilde{\rho} (\lambda)$, defined by relation
\[
  d \tilde{\rho} (\lambda) = d \rho (\lambda)  + \sum_{k=1}^n (b_{z_k}^{-1} - a_{z_k}^{-1}) \delta(\lambda - \lambda_{z_k}) d \lambda
\]
also is spectral. Moreover, there exists unique self-adjoint canonical Dirac operator $\tilde{L}$ generated on half axis by the differential
expression $\tilde{l} = B \cfrac{d}{dx} + \tilde{\Omega}(x)$ and the boundary condition \eqref{eq2.2}, for which $\tilde{\rho} (\lambda)$
is spectral function.
Wherein, the potential function $\tilde{\Omega}(x)$ is

\[
\tilde{\Omega}(x) =   \left(
                        \begin{array}{cc}
                          p(x, n) & q(x, n) \\
                          q(x, n) & -p(x, n) \\
                        \end{array}
                      \right),
\]
where $ p(x, n)$ and  $q(x, n)$ are defined by the following formulae:
\[
\begin{array}{c}
  p(x, n)  =  - \cfrac{1}{\det S(x, n)} \displaystyle \sum_{k=1}^n \displaystyle \sum_{p=1}^2 V_{z_k,(3-p)}(x) \det S_{p}^{(k)}(x, n) ,\\
  \\
  q(x, n)  = x + \cfrac{1}{\det S(x, n)} \displaystyle \sum_{k=1}^n \displaystyle \sum_{p=1}^2 (-1)^{p-1} V_{z_k,p}(x) \det S_{p}^{(k)}(x, n) ,
\end{array}
\]
where $S(x, n)$ is $n \times n$ square matrix $ S(x, n) = \{ \delta_{z_i z_j} + ( b_{z_i}^{-1} - a_{z_i}^{-1}) v_{z_i z_j} (x) \}_{i,j=1}^n $
and $S_{p}^{(k)}(x, n)$ are matrices, which are obtained from the matrix $S(x, n)$, when we replace $k$-th column of $S(x, n)$
by $H_p(x, n) = \{ - ( b_{z_j}^{-1} - a_{z_j}^{-1}) V_{z_i,p}(x) \}_{i = 1}^n$ column, $p = 1, 2$.
And for the eigenfunctions $\tilde V_m(x)$ ($ m \in \mathbb{Z}$) we obtain the representations
\[
  \tilde V_m(x) = \left(
                  \begin{array}{c}
                    V_{m,1}(x) + \cfrac{1}{\det S(x,n)} \displaystyle \sum_{k=1}^n v_{z_k m}(x) \det S_{1}^{(k)}(x, n)  \\
                    \\
                    V_{m,2}(x) + \cfrac{1}{\det S(x,n)} \displaystyle \sum_{k=1}^n v_{z_k m}(x) \det S_{2}^{(k)}(x, n)  \\
                  \end{array}
                \right).
\]

\end{theorem}

Now we want to add any finite number of new real eigenvalues $\mu_k \neq \lambda_m, \ m \in \mathbb{Z}$, to the spectra, with positive norming constants $c_k, \ k = 1, 2, \ldots, n$.

\begin{theorem}\label{thm3.4}
Let $\rho (\lambda)$ is the spectral function of the operator $L$, then the function $\tilde{\rho} (\lambda)$, defined by relation
\[
  d \tilde{\rho} (\lambda) = d \rho (\lambda)  + \sum_{k=1}^n c_{k}^{-1} \delta(\lambda - \mu_{k}) d \lambda
\]
also is spectral. Moreover, there exists unique self-adjoint canonical Dirac operator $\tilde{L}$ generated on half axis by the differential
expression $\tilde{l} = B \cfrac{d}{dx} + \tilde{\Omega}(x)$ and the boundary condition \eqref{eq2.2}, for which $\tilde{\rho} (\lambda)$
is spectral function.
Wherein, the potential function $\tilde{\Omega}(x)$ is

\[
\tilde{\Omega}(x) =   \left(
                        \begin{array}{cc}
                          p(x, n) & q(x, n) \\
                          q(x, n) & -p(x, n) \\
                        \end{array}
                      \right),
\]
where $ p(x, n)$ and  $q(x, n)$ are defined by the following formulae:
\[
\begin{array}{c}
  p(x, n)  =  - \cfrac{1}{\det S(x, n)} \displaystyle \sum_{k=1}^n \displaystyle \sum_{p=1}^2 W_{k,(3-p)}(x) \det S_{p}^{(k)}(x, n) ,\\
  \\
  q(x, n)  = x + \cfrac{1}{\det S(x, n)} \displaystyle \sum_{k=1}^n \displaystyle \sum_{p=1}^2 (-1)^{p-1} W_{k,p}(x) \det S_{p}^{(k)}(x, n) ,
\end{array}
\]
and where $W_{k}(x) := \psi(x, \mu_{k}, 0, \Omega_0), \ k = 1, 2, \ldots, n, $ and
$S(x, n)$ is $n \times n$ square matrix $ S(x, n) = \{ \delta_{i j} + c_{j}^{-1} w_{i j}(x) \}_{i,j=1}^n $ ($w_{i j}(x) = \int_{0}^{x} W^{*}_i(s) W_j(s) ds $),
and $S_{p}^{(k)}(x, n)$ are matrices, which are obtained from the matrix $S(x, n)$, when we replace $k$-th column of $S(x, n)$
by $H_p(x, n) = \{ - c_{i}^{-1} W_{i,p}(x) \}_{i = 1}^n$ column, $p = 1, 2$.
For the eigenfunctions $\tilde V_m(x)$ (for $m \in \mathbb{Z}$ we obtain the representations
\[
  \tilde V_m(x) = \left(
                  \begin{array}{c}
                    V_{m,1}(x) + \cfrac{1}{\det S(x,n)} \displaystyle \sum_{k=1}^n \int_{0}^{x} W^{*}_k(s) V_m(s) ds \det S_{1}^{(k)}(x, n)  \\
                    \\
                    V_{m,2}(x) + \cfrac{1}{\det S(x,n)} \displaystyle \sum_{k=1}^n \int_{0}^{x} W^{*}_k(s) V_m(s) ds \det S_{2}^{(k)}(x, n)  \\
                  \end{array}
                \right),
\]
and for the eigenfunctions $\tilde W_k(x)$ (for $k = 1, 2, \ldots, n$) we obtain the representations

\[
  \tilde W_l(x) = \left(
                  \begin{array}{c}
                    W_{l,1}(x) + \cfrac{1}{\det S(x,n)} \displaystyle \sum_{k=1}^n w_{k l}(x) \det S_{1}^{(k)}(x, n)  \\
                    \\
                    W_{l,2}(x) + \cfrac{1}{\det S(x,n)} \displaystyle \sum_{k=1}^n w_{k l}(x) \det S_{2}^{(k)}(x, n)  \\
                  \end{array}
                \right).
\]

\end{theorem}

The proofs of the theorems \eqref{thm3.3} and \eqref{thm3.4} are similar to the proof of the theorem \eqref{thm3.2}.
Thus, we have proven, that one can perturb the linear potential of canonical Dirac operator by adding, subtracting finite number of the eigenvalues and/or changing finite number of norming constants with having changed potential function in explicit form.

\

\noindent {\bf Remark.}
We take the operator $L(0, x, 0)$ as a "model" operator for perturbing spectral function.
Analogues theorems can be proven for the second model operator $L(0, x, \pi / 2)$.

\

\noindent {\bf Acknowledgement.} This work was supported by the RA MES State Committee of Science, in the frames of the research project № 15T-1A392.


\begin{thebibliography}{19}


\bibitem{GaLe} Gasymov MG, Levitan BM.
An inverse problem for Dirac system.
DAN SSSR, vol. 167, N. 5; 967-970pp., 1966, (in Russian).

\bibitem{Ma1} Marchenko VA.
Sturm-Liouville operators and its applications.
 Naukova Dumka, Kiiv; 1977, (in Russian).

\bibitem{Na} Naimark MA.
Linear Differential Operators.
Nauka, Moscow; 1969, (in Russian).

\bibitem{Mar} Martinov VV.
Direct methods of qualitative spectral analysis for first order non self-adjoint systems of differential equations.
Differential equations, vol. 4, N8; 1494-1508 pp., and N12; 2243-2257, 1968, (in Russian).

\bibitem{LeSa} Levitan BM, Sargsyan IS.
Introduction to spectral theory
Moscow, Nauka; 1970, (in Russian).

\bibitem{Ha2} Harutyunyan T.N.
The Cauchy problem for canonical Dirac system.
Uchenie zapiski of Artsakh University, vol. 1 (8), 1-5pp; 2004, (in Russian).

\bibitem{Sch} Schwartz L.
M$\acute{e}$thodes Mathematiques pour les sciences physiques.
Hermann, 115, Paris VI; 1961.

\bibitem{Ha3} Harutyunyan T.N.
The canonical Dirac operator with a partially given spectrum.
 Erevan. Gos. Univ. Uchen. Zap. Estestv. Nauki, N1(161), 11-19pp., 1986.

\bibitem{Ha} Harutyunyan T.N.
The transformation operators for canonical Dirac systems.
Differentialnie uravneniya, vol. 44, N8, 1011-1021 pp., 2008, (in Russian).
English translation. Differential equations, vol 44, N8, 1-12 pp., 2008.

\bibitem{AlHrMy} Albeverio S, Hryniv R, Mykytyuk Ya.
Inverse spectral problems for Dirac operators with summable potential.
Russian Journal of Math Physics, vol.12, N5; 406-423, 2005.
And with the same title in ArXiv.org; 1-25pp., February 2, 2008.



\end{thebibliography}
\end{document}